\documentclass[10pt]{article}

\usepackage{latexsym,color,amsmath,amsthm,amssymb,amscd,amsfonts}

\newtheorem{theorem}{Theorem}
\newtheorem{corollary}{Corollary}

\newtheorem{lemma}{Lemma}

\newtheorem{proposition}{Proposition}
\newtheorem{definition}{Definition}
\newtheorem{example}{Example}

\def\R{{\mathbb R}}
\def\K{{{\mathcal K}_n}}

\def\inte{{\rm int}}
\def\eps{{\varepsilon}}

\def\conv{{\rm conv}}
\def\vol{{\rm vol}}

\begin{document}

\title{Affine invariant points.
}

\author{Mathieu Meyer, Carsten Sch\"utt  and Elisabeth M. Werner 
\thanks{Partially supported by an NSF grant, a FRG-NSF grant and  a BSF grant}}

\date{}

 \maketitle

\begin{abstract}
We answer in the negative  a question by Gr\"unbaum who asked if there exists a finite basis of affine invariant points.
We give a positive  answer 
to  another question by  Gr\"unbaum about the ``size" of the set of all 
affine invariant points. Related,  we show that the set of all convex bodies $K$,  for which the set of affine invariant points is all of $\mathbb{R}^n$,
is dense in the set of convex bodies. Crucial  to establish these results, are new affine  invariant points,  not previously considered in the literature. 
\end{abstract}
\vskip 4mm
\section {Introduction.}
A number of highly influential works (see, e.g., \cite{Ga3, GaZ, GrZh},  
\cite{Hab}-\cite{Klain2}, \cite{Lud2}-\cite{LYZ2004}, 
\cite{Paouris2006, Schuster2010, 
SA1, SA2, Z3})
has directed much of the research in the theory of convex bodies  to the study of the affine geometry of these bodies. Even  questions that had been considered Euclidean in nature,  turned out to be affine problems - among them the famous Busemann-Petty Problem (finally laid to rest in
\cite{Ga1, GaKoSch,  Z1, Z2}).
\par
The affine structure of convex bodies is closely related to the symmetry structure of the bodies.
From an affine point of view, ellipsoids  are the most symmetric convex bodies,  and simplices are considered to be among the  least symmetric ones.  This is reflected in  many affine invariant inequalities (we give examples below)  where ellipsoids and simplices
are the  extremal cases.
However, simplices have many affine symmetries. Therefore, a more  systematic study
for  symmetry of convex bodies is needed. Gr\"unbaum,  in his seminal paper \cite{Gruenbaum1963}, initiated such a study. A crucial notion in his work,  the {\em affine invariant point},  allows to  analyze the symmetry situation. In a nutshell: the more affine invariant points, the fewer symmetries.
\par
In this paper, we address several issues that were left open in 
Gr\"unbaum's paper. For instance, it was not even 
known 
whether there are ``enough" affine invariant points. We settle  this  in Theorem \ref{dicht} below. 
\par
Let $\mathcal K_n$ be  the set of all convex bodies in 
$\mathbb R^{n}$ (i.e., compact convex subsets of $\mathbb R^{n}$
with nonempty interior).  Then (see Section 2 for the  precise definition)
a  map 
$p:{\mathcal K}_n \rightarrow\mathbb R^{n}$ is called an affine invariant point,  if $p$ is continuous and if
for every nonsingular affine map $T:\mathbb R^{n}\rightarrow \mathbb R^{n}$ one has, 
\begin{equation*}
p(T(K))=T(p(K)).
\end{equation*}
An important example of an affine invariant point is the centroid $g$. More examples will be given throughout the paper.
Let  $\mathfrak{P}_n$ be  the set of affine invariant points on ${\mathcal K}_n$,
\begin{equation*} 
\mathfrak{P}_n =\{ p: \mathcal K_n \rightarrow\mathbb R^{n}  \big|  \  p \  \text{ is  affine invariant} \}.
\end{equation*}
Observe  that $\mathfrak{P}_n$ is an affine subspace of $C(\mathcal{K}_n, \mathbb{R}^n)$, the continuous functions on $\mathcal{K}_n$ with values in $  \mathbb{R}^n$.  We denote by $V\mathfrak{P}_n$ the subspace parallel to $\mathfrak{P}_n$. Thus, with  the centroid $g$, 
\begin{equation*}\label{vektorraum}
V\mathfrak{P}_n = \mathfrak{P}_n - g.
\end{equation*}
Gr\"unbaum \cite{Gruenbaum1963}  posed the problem if  there  is a finite basis of affine invariant points, i.e. affine invariant points 
$p_i\in \mathfrak{P}_n$, $1 \leq i \leq l$,  such that every  $p\in \mathfrak{P}_n$ can be written as $$p=\sum_{i=1}^l \alpha_i p_i, \hskip 3mm\text{ with } \alpha_i \in \mathbb{R}\hskip 2mm\text{ and  }\sum_{i=1}^l \alpha_i =1.$$
We answer this question in the negative  and prove:
\begin{theorem}\label{unendlich}
$V\mathfrak{P}_{n}$ is infinite dimensional for all $n\geq2$.
\end{theorem}
\par
\noindent
In fact, we will see that, with a suitable norm,  $V\mathfrak{P}_n $ is a Banach space.
Hence, by Baire's theorem,  a basis of $\mathfrak{P}_{n}$ is not even countable.
\par
For a fixed body $K \in \mathcal{K}_n$, we let
$$
\mathfrak{P}_n(K)=\{p(K):   p \in \mathfrak{P}_n\}.
$$
Then Gr\"unbaum conjectured \cite{Gruenbaum1963} that
for every $K \in \mathcal{K}_n$, 
\begin{equation}\label{frage2}
\mathfrak{P}_n(K) = \mathfrak{F}_n(K),
\end{equation}
where
$\mathfrak{F}_n(K)= \{ x \in \mathbb{R}^n: Tx=x, \text{  for all affine  } T 
\text{ with  } TK=K\}$. 
We give a positive  answer to  this conjecture,  when $ \mathfrak{P}_n(K) $ is $(n-1)$-dimensional.
Note also that if $K$ has {\em enough symmetries},
in the sense that 
${\mathfrak F}_n(K)$
is reduced to one point $x_K$, then $ \mathfrak{P}_n(K)  = \{x_K\}$.
\vskip 2mm
\begin{theorem} \label{Pn=Fn}
Let $K \in \mathcal K_n$ be such that $ \mathfrak{P}_n(K) $ is (n-1)-dimensional. Then
$$
\mathfrak{P}_n(K) = \mathfrak{F}_n(K).
$$
\end{theorem}
Symmetry or enough symmetries, are key in many problems. The affine invariant inequalities connected with the
affine geometry often have ellipsoids,  respectively simplices as extremal cases. Examples  are the $L_p$ affine isoperimetric
inequalities of the $L_p$ Brunn Minkowski theory, a theory  initiated by Lutwak in the groundbreaking paper \cite{Lu2}. For related results we refer  to e.g.  \cite{BLYZ2012, MeyerWerner1998, MeyerWerner2000}, \cite{SchusterWeberndorfer2012}-\cite{SW5},  \cite{WernerYe2008}, \cite{WernerYe2010}. The corresponding $L_p$ affine isoperimetric
inequalities, established by Lutwak \cite{Lu2} for $p>1$ and in \cite{WernerYe2008} for all other $p$ - the case $p=1$ being the classical
affine isoperimetric inequality \cite{BlaschkeBook} - 
 are stronger  than the 
celebrated  {\em Blaschke Santal\'o inequality} (see e.g., \cite{GardnerBook, SchneiderBuch};
and e.g., \cite{Boroczky2010, NPRZ2010} for recent results):  the volume product  of polar reciprocal  convex bodies
is maximized precisely by ellipsoids. 
\par
It is an open problem  which convex bodies are minimizers for the  Blaschke Santal\'o inequality.
 Mahler   conjectured that the minimum is attained
for the simplex. A major  breakthrough towards Mahler's conjecture is the {\em inequality of Bourgain-Milman}
\cite{BourgainMilman1987}, which has been reproved with completely
different methods by Kuperberg \cite{Kuperberg2008} and by Nazarov \cite{ Nazarov}. See also \cite{GMR1988, R1, R2, S-R} for related results. Even more  surprising is that   it  is not   known whether the minimizer is a polytope. 
The strongest indication to date that it is indeed the case is given in \cite{ReisnerSchuettWerner2011}.
\par
Another example  is the {\em Petty projection inequality} \cite{Petty1972},  a far stronger inequality than the classical
isoperimetric inequality, and its $L_p$ analogue,  the {\em $L_p$ Petty projection inequality}, established by
Lutwak, Yang, and Zhang \cite{LutwakYangZhang200/1} (see also Campi and Gronchi \cite{CampiGronchi}). These inequalities were recently strengthened and extended by Haberl and Schuster \cite{HabSch2}. It is precisely the  ellipsoids that are  maximizers in all these inequalities. On the other hand, the reverse of the Petty projection inequality, the {\em Zhang projection inequality} \cite{Zhang1991}, has the simplices as maximizers.
\par
Gr\"unbaum \cite{Gruenbaum1963} also asked , whether $\mathfrak{P}_n(K)= \mathbb{R}^n$, if $\mathfrak{F}_n(K)= \mathbb{R}^n$.
A first step toward solving this problem, is to clarify 
if there is a convex body $K$ such that
$\mathfrak{P}_n(K)=\mathbb{R}^n$. Here, we answer this question in the affirmative and prove 
that the set of all $K$ such that $\mathfrak{P}_n(K)=\mathbb{R}^n$,  is dense in $\mathcal{K}_n$ and consequently
the set of all $K$ such that $\mathfrak{P}_n(K)=\mathfrak{F}_n(K)$,  is dense in $\mathcal{K}_n$. 

\begin{theorem} \label{dicht}
The set of all $K \in \mathcal{K}_n$ such that $\mathfrak{P}_n(K) = \mathbb{R}^n$ is open and dense in $\left(\mathcal{K}_n, d_H\right)$.
\end{theorem}
\noindent
Here, $d_H$ is the  Hausdorff metric on $\mathcal K_n$, defined as
\begin{equation}\label{Hausdorff}
d_H(K_1, K_2) = \min\{\lambda \geq 0: K_1 \subseteq K_2+\lambda B^n_2;   K_2 \subseteq K_1+\lambda B^n_2 \},
\end{equation}
where $B^n_2$ is the Euclidean unit ball centered at $0$. More generally, $B^n_2(a,r)$, is the Euclidean ball centered at $a$ with radius $r$.  
We shall use the following  well known fact. 
\newline
Let $K_m, K\in \K$. Then $d_H(K_m,K)\to 0$ if and only if for some  $\varepsilon_m\to 0$ one has 
\begin{equation}\label{Konvergenz}
(1-\eps_m) \big(K-g(K)\big)\subset K_m-g(K_m)\subset (1+\eps_m)\big(K-g(K)\big)\hbox{ for every $m$}.
\end{equation}

\par
To establish Theorems \ref{unendlich} - \ref{dicht}, we need to introduce new examples of affine invariant points,  that have not previously been considered in the literature.

\vskip 4mm
\section { Affine invariant points and sets: definition and  properties.}
Let $K \in \mathcal K_n$. Throughout the paper,  $\mbox{int}(K)$ will denote  the interior,  and 
 $\partial K$ the boundary of $K$. The $n$-dimensional volume of $K$ is $\text{vol}_n(K)$, or simply $|K|$.
$K^\circ=\{y \in \mathbb R^{n}: \langle x,y \rangle \leq 1  \   \forall x \in K\}$
 is the {\it polar body of $K$}. 
More generally, for $x$ in $\mathbb R^{n}$,   the {\it  polar of $K$ with respect to $x$} is 
$K^x=(K-x)^\circ +x$.   
\par
A map $p:\mathcal K_n \rightarrow\mathbb R^{n}$
is said to be continuous if it is continuous when $\mathcal K_n $ is
equipped  with the Hausdorff metric and $\mathbb R^{n}$ with
the Euclidean norm $\| \cdot \|$.
\par
Gr\"unbaum \cite{Gruenbaum1963} gives the following definition of affine invariant points. Please note that
formally we are considering maps, not points.

\vskip 3mm
\begin{definition}\label{aip}
A  map 
$p:{\mathcal K}_n \rightarrow\mathbb R^{n}$ is called an affine invariant point,  if $p$ is continuous and if
for every nonsingular affine map $T:\mathbb R^{n}\rightarrow \mathbb R^{n}$ one has
\begin{equation} \label{def1}
p(T(K))=T(p(K)).
\end{equation}
Let  $\mathfrak{P}_n$ the set of affine invariant points in $\mathbb{R}^n$,
\begin{equation} \label{def:aip}
\mathfrak{P}_n =\{ p: \mathcal K_n \rightarrow\mathbb R^{n}  \big|  \  p \  \text{ is  affine invariant} \}, 
\end{equation}
and for a fixed body $K \in \mathcal{K}_n$,
$
\mathfrak{P}_n(K)=\{p(K):   p \in \mathfrak{P}_n\}.
$
\vskip 2mm
We say an affine invariant point $p\in  \mathfrak{P}_n $  proper,  if for all
$K\in {\mathcal K}_n$, one has
$$
p(K)\in \operatorname{int}(K).
$$
\end{definition}
\vskip 2mm
\noindent
{\bf Examples.}
Well known examples (see e.g. \cite{Gruenbaum1963})  of  proper   affine invariant points  of a convex body
$K$ in
$\mathbb R^{n}$ are
\newline
(i) the {\em centroid} 
\begin{equation}\label{centroid}
g(K)
=\frac{\int_{K}xdx}{|K|};
\end{equation}
\newline
(ii)  the {\em Santal\'o point}, the unique point
$s(K)$
for which the volume product 
$
|K||K^{x}|
$
attains its minimum;
\newline
(iii) the center $j(K)$ of the {\em ellipsoid of maximal volume}
$\mathcal{J}(K)$ contained in $K$,  or {\em John ellipsoid} of $K$; 
\newline
(iv)
the center $l(K)$ of the {\em ellipsoid of minimal volume} $\mathcal{L}(K)$
containing $K$, or {\em L\"owner ellipsoid} of $K$. 
\vskip 2mm
Note  that if $T(K)=K$ for some affine map $T:\mathbb R^{n}\rightarrow \mathbb R^{n}$ and some
$K\in {\mathcal K}_n$, then for every $p\in \mathfrak{P}_n$, one  has $p(K)=p(T(K))=T(p(K))$. It follows that if $K$ is 
centrally symmetric or is a simplex, then $p(K)=g(K)$  for every $p\in \mathfrak{P}_n$, hence ${\mathfrak P}_n(K)=\{g(K)\}$.
\vskip 2mm
The continuity property is an essential part of  Definition \ref{aip}, as,  without it, pathological
 affine invariant points can be constructed.
The next example  illustrates this.  
\vskip 2mm
\begin{example}
\label{example2}
Let $\mathcal{P}_n$ be the set of all convex polytopes in  $\mathcal{K}_n$ and define for $ P\in \mathcal{P}_n$, 
$$p(P)=\frac{1}{m} \sum_{i=1}^m v_i(P),$$
where $v_1(P), \dots v_m(P)$ are the vertices of $P$.
For $K\in \mathcal{K}_n\setminus \mathcal{P}_n$, let
$p(K)=g(K)$, the centroid of $K$.
Then $p:  \mathcal{K}_n\to  \mathbb{R}^n$ is affine invariant, but it is not continuous at any point.
\end{example}
\par
Indeed, let $K \in \K$. We approximate $K$ by a polytope $P$, and, in turn, approximate $P$ by a polytope $P_l$
by replacing one vertex $v$ of $P$ by sufficiently  many vertices $v_1, \dots, v_l$ near $v$. When $l \rightarrow \infty$, 
$p(P_l) \rightarrow v \in \partial P$. $P_l$ is near $K$, but $p(K) =g(K)$.

\vskip 3mm
Next,  we introduce  the notion of affine invariant set mappings, or, in short, affine invariant sets.
There,  continuity of a map $A:\mathcal K_n \rightarrow \mathcal K_n$
is meant when $\mathcal K_n $ is
equipped on both sides with the Hausdorff metric. Our Definition \ref{ais} of affine invariant sets  differs  from the one given by Gr\"unbaum \cite{Gruenbaum1963}.
\par
\begin{definition}\label{ais}
A  map 
$A: \mathcal K_n \rightarrow \mathcal K_n $ is  called an affine invariant set mapping,  if $A$ is continuous and if 
 for every nonsingular affine map $T$ of $\mathbb R^{n}$, one has
\begin{equation*} \label{def2}
A(TK)=T(A(K)).
\end{equation*}
We then  call $A(K)$,  or  simply the map $A$, an affine invariant set mappings.
We denote by $\mathfrak{S}_n$ the set of affine invariant set mappings, 
\begin{equation} \label{def:ais}
\mathfrak{S}_n =\{ A: \mathcal K_n \rightarrow \mathcal K_n \big| A \  \text{ is  affine invariant and continuous} \}.
\end{equation}
\vskip 2mm  
We say that $A\in \mathfrak{S}_n$ is proper,  if 
$A(K)\subset \inte(K)$ for every $K\in  \mathcal K_n$.
\end{definition}
\par
Known examples (see e.g.,  \cite{Gruenbaum1963})  of  affine invariant sets are 
the  John ellipsoid and  the  L\"owner ellipsoid.
Further examples will be given all along this paper.
\vskip 3mm
\noindent
{\bf Remarks.} 
(i) It is easy to see  that if $ \lambda \in \mathbb{R}$, $p, q  \in \mathfrak{P}_n$ and $A  \in \mathfrak{S}_n$, then
$p\circ A  \in {\mathfrak P}_n$ and $ (1-\lambda)p+\lambda q\in \mathfrak{P}_n$. 
Thus,  $\mathfrak{P}_n$ is an affine space and for every $K \in \mathcal{K}_n$,  ${\mathfrak P}_n(K)$ is an affine subspace of $\mathbb{R}^n$.
Moreover, for  $A, B  \in \mathfrak{S}_n$,  the maps
$$K \rightarrow (A \circ B)(K),  \hskip 3mm (1-\lambda)A(K) + \lambda B(K)    
\hskip 3mm \text{and}    \hskip 3mm \conv[A ,B] (K)=\conv[A(K), B(K)]$$
are  affine invariant set mappings.
\par
(ii)
Properties (\ref{def1}) and (\ref{def2}) imply in particular that  for every translation by a fixed vector $x_0$ and  for every
convex body $K \in \mathcal{K}_n$, 
\begin{equation} \label{trans}
p(K + x_0) = p(K) + x_0, \  \text{for every } \  p \in {\mathfrak P}_n
\end{equation}
and 
\begin{equation} \label{trans1}
A(K + x_0) = A(K) + x_0,  \  \text{for every } \  A \in {\mathfrak S}_n
\end{equation}
\par
(iiii) Unless $p=q$, it is not possible to compare two different affine invariant points $p$ and $q$ via an inequality of the following type
\begin{equation} \label{compare}
\| p(K) - p(L) \| \geq c\  \| q(K) - q(L) \|,
\end{equation}  
where $\| \cdot\|$ is a norm on $\mathbb{R} ^n$ and $c>0$ a constant. Indeed, 
by (ii),  $p(K - p(K)) = 0$ and $q(L - q(L)) = 0$. Therefore, if  (\ref{compare}) would hold, then
\begin{eqnarray*} 
0 &=& \| p(K - p(K)) - p(L - p(L)) \| \geq c\  \| q(K -p(K)) - q(L - p(L)) \| \\
&=&c\   \| q(K) - p(K) -q(L) + p(L)\|.
\end{eqnarray*}
Choose  now for $L$ a symmetric convex body. Then   $ \| q(K) - p(K) \| =0$, or $p(K)=q(K)$.

\vskip 2mm
Remark (i) provides  examples of non-proper affine invariant points:
once there are two different affine
invariant points, there are affine invariant points  $p(K)\notin K$, i.e.  non-proper affine invariant points.
An explicit example  is the convex body $C_n$ constructed in \cite{MeyerSchuettWerner2011}, for which  the centroid  and the Santal\'o point   differ.
\vskip 1mm
\noindent
The next results describe some properties   of affine invariant points and sets.
\vskip 2mm
\noindent
\begin{proposition} \label{lemma2bis}
Let $p,q\in  {\mathfrak P}_n$ and suppose that  $p$ is proper. 
For $K\in \mathcal K_n$, define
$$\phi_q(K)=\inf \left\{t\ge 0: q(K)-p(K)\in t\left(K-p(K)\right)\right\}.$$
Then  $\phi_q:  \mathcal K_n\to {\mathbb R}_+$ is continuous and
\par
(i) there exist $c=c(q) > 0$ such that 
$$q(K)-p(K)\in c\ \left(K-p(K)\right) \hbox{ for every } K\in \mathcal K_n. $$
\par
(ii) If moreover $q$ is proper, then one can chose $c \in (0,1)$ in (i).
\end{proposition}
\noindent
{\bf Proof.}
Since $p(K)\in \inte(K)$, 
${\mathbb R}^n=\cup_{t\ge 0} t(K-p(K))$. Therefore $\phi_q$ is well defined. 
Now we show that $\phi_q$ is continuous.
Suppose that $K_m\to K$ in $( \mathcal K_n, d_H)$. By definition, we have
$$q(K_m)-p(K_m)\in \phi_q(K_m) \left(K_m-p(K_m)\right)\hbox{  for all $m$}.$$
By continuity of $p$ and $q$ it follows that 
$$q(K)-p(K)\in \liminf_m \phi_q(K_m)\ \left( K-p(K)\right),$$
and thus $$\phi_q(K)\le  \liminf_m\phi_q(K_m).$$
Since $p$ is proper,  there exists $d>0$, such that $B_2^ n\subseteq
d(K-p(K))$. Since $K$ is bounded, there exists  $D>0$ such that 
$d(K-p(K))\subseteq DB_2^ n$. Let $\eta >0$ and fix $\eps=\eta/d>0$. Since $K_m\to K$,  and by  continuity of $p$ and $q$, there exists $m_0>0$ such that for every
$m\ge m_0$, 
$$K-p(K)\subseteq K_m-p(K_m)+\eps B_2^n\subseteq K_m-p(K_m)+\eta(K-p(K))$$
and $$q(K_m)-p(K_m)\in q(K)-p(K)+ \eps B_2^n\subseteq (\phi_q(K)+\eta)(K-p(K)).$$
Now we observe that,  if two convex bodies $A$ and $ B$ in  ${\mathbb R}^n$ satisfy $A\subseteq B+tA$
for some $0<t<1$, then $A\subseteq B/(1-t)$.  
It then follows  that for every
$m\ge m_0$, 
$$K-p(K)\subseteq \frac{1}{1-\eta} (K_m-p(K_m)).$$
Hence
$$q(K_m)-p(K_m)\in \frac{ \phi_q(K)+\eta}{1-\eta}\left(K_m-p(K_m)\right), $$
and thus $$\phi_q(K)\ge \limsup_m\phi_q(K_m),$$
and the continuity of $\phi_q$ is proved. Assertions
(i) and (ii)  follow from the continuity of $\phi_q$. Indeed,  by affine invariance, we may reduce the problem to the set $\{ K\in \mathcal K_n: B_2^n\subseteq K\subseteq nB_2^n\}$, which is compact in   $\mathcal K_n$. $\Box$
\vskip 3mm
\begin{lemma}
Let $p, q\in {\mathfrak P}_n$  and suppose that $p$  is proper. Then there exists a proper $r\in {\mathfrak P}_n$ such that $q$ is an affine combination of $p$ and $r$.
\end{lemma}
\noindent
{\bf Proof.} By the preceding proposition, there is  $c>0$ such that 
$q(K)-p(K)\in c (K-p(K)$ for all $K \in  \mathcal K_n$.
Put $r= \frac{q-p}{2c}+p$. Then $r \in {\mathfrak P}_n$ and  $
q=2c r+ (1-2c)p$ is an affine combination of $p$ and $r$. Since $p(K)\in \inte(K)$,
$$r(K)\in \frac{1}{2} (K+p(K)) \subseteq \inte(K),  \  \  \text{for all}  \  \   K \in \mathcal K_n.\hskip 10mm \square $$

\vskip 3mm
Analogous results to Proposition \ref{lemma2bis} for  affine invariant sets are also valid. We omit their proofs.
\vskip 2mm
\begin{proposition} \label{afflemma}
Let $A\in \mathfrak{S}_n$, 
$p, q \in \mathfrak{P}_n$ and suppose that $p$ is proper.
Then there exists a constant $c_1>0$ such that  for every $ K\in\mathcal{K}_n$, 
$$A(K)- q\big(K\big)\subseteq c_1\  (K-p(K)).$$
If moreover $A$ is proper and $p=q$, one can choose $c_1<1$.
\end{proposition}

\vskip 3mm

\begin{lemma} Let $A\in\mathfrak{S}_n$ and $p\in\mathfrak{P}_n$ be  proper.  Then there exists $t>0$ such that 
$$K\to t(A(K)-p(K))+ p(K)=tA(K)+(1-t)p(K)$$
is a proper affine invariant set mapping.
\end{lemma}

\vskip 3mm
The next proposition  gives a reverse inclusion for affine invariants sets. We need  first another lemma, where, as in the proposition, $g$ denotes the center of gravity.
\vskip 3mm
\begin{lemma}\label{boule} For every $D,d >0$ and $n\ge 1$, there exists $c>0$ such that,  whenever $K\in \K$ satisfies
$K\subseteq DB_2^n$ and $|K|\ge d$, then $cB_2^n\subseteq K-g(K)$.
\end{lemma}
\noindent
{\bf Proof.} Suppose that $K\in \K$ satisfies the two assumptions. Define
$$ c_K=\sup\{c \geq 0: \ c B_2^n\subseteq K-g(K)\}.$$ Then
$c_K>0$, and there exists $x\in \partial K$ such that $\|x-g(K)\|=c_K$. Since $K-g(K)\subseteq n\big(g(K)-K\big)$,
the length of the chord of $K$ passing through $g(K)$ and $x$ is not bigger than $(n+1)c_K$. 
Let $u\in S^{n-1}$ be the direction of the segment $[g(K), x]$ and let $P_u K$ be the orthogonal projection of $K$ onto $u^{\perp}$,
the subspace orthogonal to $u$. Then, 
$$d\le |K|\le \|(g +\R u)\cap K\|\ |P_u K| \le (n+1) c_K \ D^{n-1} |B_2^{n-1}|.$$
The second inequality follows from a result by Spingarn \cite{spingarn}.
Thus we get a strictly positive  lower bound $c$ for $c_K$ which depends only on $n$, $d$ and $D$. $\square$

\vskip 3mm

\begin{proposition} Let $A$ be an affine invariant set mapping. Then there exist $c>0$ such that
$$c\big(K-g(K)\big)\subseteq A(K)-g\big(A(K)\big),  \   \hbox{ for every } K\in \K.$$
\end{proposition}
\vskip 2mm
\noindent
{\bf Proof.} We first prove that there exists $d>0$ such that 
$|A(K)|\ge d|K|$  for every $K\in \K$.
By affine invariance, it is enough  to prove that 
$$\inf_{\{K\in \K: B_2^n \subset K\subseteq nB_2^n\}} \frac{|A(K)|}{|K|} >0.$$
Since $K\to \frac{|A(K)|}{|K|}$ is continuous and since $\{K\in \K: B_2^n \subseteq K\subset nB_2^n\}$ is compact in $\K$,
this infimum  is  a  minimum  and it is strictly positive. 
By  Proposition \ref{afflemma}, applied with $q=g \circ A$, there exists $c>0$ such that 
$$A(K)-g\big(A(K)\big)\subseteq c\big( K-g(K)\big), \  \   \hbox{ for every } K\in \K.$$
Therefore, 
$$A(K)-g\big(A(K)\big)\subseteq 2ncB_2^n, \  \  \hbox{ for every }K\in {\mathcal K}.$$
By Lemma \ref{boule}, there is  $c_0>0$ such that 
$$c_0B_2^n\subseteq A(K)-g\big(A(K)\big), \  \  \hbox{ for every }K\in {\mathcal K}.$$
Now, since $K-g(K)\subseteq 2n B_2^n $ for every $K\in  \K$ with $ B_2^n \subseteq K\subset nB_2^n$, we get the result for $K\in  \K$ with $ B_2^n \subseteq K\subset nB_2^n$, and thus
for all $K\in \K$ by affine invariance. $\square$

\section{Several  questions by Gr\"unbaum.}
We now give the proof of Theorems \ref{unendlich} - \ref{dicht}. To do so, we first need to introduce new affine invariant points.

 \vskip 3mm 
\subsection {The convex floating body as an affine invariant set mapping.}

Let $K\in{\mathcal K}_n$ and $0 \leq \delta  <\left(\frac{n}{n+1}\right)^n$.  For $u \in \mathbb{R}^n$ and $a \in \mathbb{R}$, $H=\{x \in \mathbb{R}^n: \langle x, u\rangle  = a\}$ is the hyperplane 
orthogonal to $u$ and $H^+=\{x \in \mathbb{R}^n: \langle x, u \rangle \geq  a\}$ and $H^-=\{x \in \mathbb{R}^n: \langle x, u \rangle  \leq  a\}$ are the two half spaces determined by $H$. 
Then the (convex) floating body 
$K_{\delta}$  
\cite{SchuettWerner1990} of $K$ is
the intersection of all halfspaces $H^+$ whose
defining hyperplanes $H$ cut off a set of volume at most $\delta |K|$
from $K$, 
\begin{equation} \label{schwimm}
K_{\delta}=\bigcap_{\{H: |H^-\cap K| \leq \delta|K|\}} {H^+}. 
\end{equation}
Clearly, $K_0=K$ and $K_{\delta}\subseteq K$ for all $\delta \geq 0$.
The condition $\delta  <\left(\frac{n}{n+1}\right)^n$ insures that $g(K)\in \inte(K_{\delta})\not=\emptyset$ (see \cite{SchuettWerner1990}).
Moreover, for
all invertible affine maps $T$, one has 
\begin{equation}\label{Affine:map:floating:body}
(T(K))_{\delta}=T\left(K_\delta \right).
\end{equation}
To prove that $K\to K_{\delta}$ is continuous from ${\mathcal K}_n$ to ${\mathcal K}_n$,  we need some notation. For $u \in S^{n-1}$, we define
$a_{\delta,K}(u)$ to be
unique real number such that 
\begin{equation}\label{adelta}
\vol_n\left(\{x\in K: \langle x,u\rangle \ge a_{\delta,K}(u)\}\right) =\delta \ \vol_n(K).
\end{equation}
Then
$$K_{\delta}=\bigcap_{u\in  S^{n-1}} \left\{x\in K: \langle x,u\rangle \le a_{\delta,K}(u)\right\}.$$

\vskip 3mm
\begin{lemma}\label{lemma:f(a)}
Let $K \in{\mathcal K}_n$, $u\in S^{n-1}$, $0<\delta<\left(\frac{n}{n+1}\right)^n$ and $f(t)= \left|\{x\in K: \langle x,u\rangle =t\}\right|$. Let
$a\in \mathbb{R}$ satisfy $\int_a^{+\infty} f(t) dt =\delta \int_{-\infty}^{+\infty} f(t) dt$. Then one has 
$$f(a)\ge \delta^{\frac{n-1}{n}} \max_{t\in \mathbb{R}} f(t).$$
\end{lemma}
\vskip 2mm
\noindent 
{\bf Proof.}  By the  Brunn-Minkowski theorem (see \cite{GardnerBook, SchneiderBuch}), $ f^{\frac {1}{n-1}}$ is concave on $\{f>0\}$. Put $M= f(m) =\max_{t\in \mathbb{R}} f(t)$. 
\par
We  suppose first that $m<a$. Let $g$ be the affine function on $\mathbb{R}$ such that $g^{n-1}(m)=f(m)$ and $g^{n-1}(a)=f(a)$. As
$ f^{\frac {1}{n-1}}$ is concave on $\{f\not=0\}$, one has $g^{n-1}\le f$ on $[m,a]$ and $g^{n-1}\ge f$ on $\{f\not=0\}\setminus [m,a]$. Thus  there exists $c\le m$ and $d\ge a$, such that $g(c)>0$, $g(d)>0$, 
$$\int_c^d g^{n-1}(t) dt =\int_{-\infty}^{+\infty} f(t) dt\hbox{\hskip 4mm and }  \hskip 4mm \int_a^d g^{n-1}(t) dt =\int_a^{+\infty} f(t) dt.$$ 
Let $g_1= g{\bf 1}_{[c,d]}$. Since $g$ is non increasing on $[c,d]$,   $g^{n-1}_1(c)\geq M=f(m)$. Moreover,  by construction, 
$g_1^{n-1}(a)=f(a)$. We replace now $g_1$ with a new function $g_2$ that is affine on its support $[c',d']$, $c'\le a\le d'$, and satisfies
$g_2^{n-1}(a)=f(a)$, $g_2^{n-1}(d')=0$,
  $$\int_{c'}^{d'}g_2^{n-1}(t) dt =\int_{-\infty}^{+\infty} f(t) dt \hbox{\hskip 4mm and } \hskip 4mm
\int_{a}^{d'}g_2^{n-1}(t) dt =\int_a^{+\infty} f(t) dt.$$ One  still  has $g_2^{n-1}(a)=f(a)$ and clearly
$g_2^{n-1}(c')\ge g_1^{n-1}(c)\ge M$.
Now, an easy computation gives 
$$f(a)= g_2^{n-1}(a)=\delta^{\frac{n-1}{n}}, \  \    g_2^{n-1}(c')\ge \delta^{\frac{n-1}{n}} M.$$
\par
 We suppose next  that $m\ge a$. The same reasoning,  with $1-\delta$ instead 
of $\delta$, gives $f(a)\ge (1-\delta)^{\frac{n-1}{n}}M$.
Since $0<\delta<\frac{1}{2}$, the statement follows. $\square$

\vskip 3mm
\noindent
\begin{proposition}
Let $0<r\leq R < \infty$ and let $K\in {\mathcal K}_n$ satisfy,    $rB_2^n\subseteq K\subseteq RB_2^n$.  Let $0<\delta <\frac{1}{2}$ and $\eta >0$. There there exists $\eps>0$ (depending only on $r, R, n,\delta$) such that,  whenever a convex body 
$L$ satisfies $d_H(K,L)\le \eps$, one has for every $u\in S^{n-1}$
$$ a_{\delta,K}(u)-\eta\le a_{\delta,L}(u)\le a_{\delta,K}(u)+\eta. $$
\end{proposition}
\vskip 2mm
\noindent 
{\bf Proof.}
Let $\rho>0$.
With the hypothesis on $K$, we may choose $\eps>0$ small enough such that whenever $d_H(K,L)\le \eps$, then $  (1-\rho)K\subseteq L\subseteq (1+\rho)K$.
Fix $u\in S^{n-1}$ and define $$ f_K(t)= \vol_{n-1}\left(\{x\in K: \langle x,u\rangle =t \} \right)\hbox{ and } f_L(t)= \vol_{n-1}\left(\{x\in L: \langle x,u\rangle =t\}\right).$$ 
Then $ a_{\delta,K} : = a_{\delta,K}(u) $ and $
a_{\delta,L} : =a_{\delta,L}(u)$ satisfy
$$\int_{ a_{\delta,K}}^{+\infty}f_K(t) dt=\delta |K| \hskip 4mm \hbox{ and } \hskip 4mm \int_{ a_{\delta,L}}^{+\infty}f _L(t) dt=\delta |L|.$$
Let $\theta>0$. For $\rho>0$ small enough one has, 
 $$|K\Delta L| \le \left((1+\rho)^n - (1-\rho)^n \right) 
 |K|\le \theta.$$ 
 For such a $\rho$ one has also 
 \begin{eqnarray*}
 \int_{\mathbb{R}} |f_K(t)-f_L(t)| dt &\le& \int_{\mathbb{R}} \vol_{n-1}\left(\{x\in K\Delta L: \langle x,u\rangle =t\}\right) dt= |K\Delta L|\le 
 \theta,
 \end{eqnarray*}
 so that 
 \begin{eqnarray*}
\left|\int_{[a_{\delta,K},a_{\delta,L}]} f_K(t) dt \right|
\le \left|\int_{a_{\delta,K}}^{+\infty} f_K(t) dt -\int_{a_{\delta,L}}^{+\infty} f_L(t) dt\right| + \int_{a_{\delta,L} }^{+\infty} |f_K(t)-f_L(t)| dt\le 2\theta.
 \end{eqnarray*}
For  $\alpha>0$ given, let $\theta= \frac{\alpha |K|}{2}$. Then 
$$\left|\int_{ a_{\delta,L}}^{+\infty}  f_K(t) dt-\delta |K| \right|\le   \alpha \ \delta |K|. $$
For some $\beta\in [-\alpha, \alpha]$, one has hence that
$a_{\delta,L}= a_{(1+\beta) \delta,K}$.
Concavity of $f_K^{\frac{1}{n-1}}$ on $\{f_K\not=0\}$ implies that
$$\left|\int_{[a_{\delta,K},a_{\delta,L}]} f_K(t) dt\right|
\ge |a_{\delta,K}-a_{\delta(1+\beta),K}| \min \left(f_K(a_{\delta,K}), 
f_K(a_{(1+\beta)\delta, K})\right).$$
If $M=\max_{t \in \mathbb{R}} f_K(t)$, we get by  Lemma \ref{lemma:f(a)},
$$\min \bigg(f_K(a_{\delta,K}), 
f_K(a_{(1+\alpha)\delta }\bigg)\ge \bigg(\min(1+\beta, 1)\delta\bigg)^{\frac{n-1}{n}}\ge \big((1-\alpha)\delta\big)^{\frac{n-1}{n}} M.$$
Since $K\subseteq R B_2^n$, we estimate $M$ from above by 
$$M\le \gamma_n =R^{n-1} |B_2^{n-1}|,$$ 
which is an upper bound independent of $u$.
It follows that if $\alpha>0$ is small enough, then
$$|a_{\delta,K}-a_{\delta,L}|\le 2 \  \theta \  \gamma_n \  \big((1-\alpha) \delta\big)^{-\frac{n-1}{n}} \le \eta. $$   
$\square$ 
\vskip 3mm 
The next proposition shows that the map $K \mapsto K_\delta$  as defined in (\ref{schwimm}), is an affine invariant  set mapping.
\vskip 3mm 
\noindent
\begin{proposition} \label{cont:Kdelta} 
For  $0<\delta<\big(\frac{n}{n+1}\big)^n$, the mapping $K\mapsto K_{\delta}$ is is an affine invariant  set  mapping from ${\mathcal K}_n$ to 
${\mathcal K}_n$ .  
\end{proposition}
\vskip 2mm
\noindent
{\bf Proof.}  
We take  $0<\delta<\big(\frac{n}{n+1}\big)^n$ so that ${\rm int} (K_\delta) \neq \emptyset$ and $g(K) \in {\rm int} (K_\delta)$ .  It is clear that 
$K \rightarrow K_\delta$ is an affine invariant  mapping and it is clear that $g(K) \in K_\delta$. 
We now fix a body $K\in {\cal K}_n$ and we verify the continuity of the  mapping  $K \rightarrow K_\delta$ at $K$. We may suppose that $0$ is the center of mass of $K$. For some $0<r\leq R < \infty$, one has
$$rB_2^n\subseteq K_{\delta} \subseteq K\subseteq RB_2^n.$$ 
By the choice of $\delta$,  $a_{\delta,K}(u)>0$ 
for every $u\in S^{n-1}$, where  $a_{\delta,K}(u) $ is as in (\ref{adelta}). 
Let $\eta, \eta'>0$ satisfy  $\eta'\le \eta r\le \eta\ \min_{u \in S^{n-1}} a_{\delta, K}(u)$. We use the notation of the preceding proposition to find $\eps>0$ such that for any $L$ with
$d_H(K,L)\le \eps$, one has   $$a_{\delta,K}(u) -\eta'\le a_{\delta,L}(u) \le a_{\delta,K}(u) +\eta',$$  
or
$$(1-\eta) a_{\delta,K}(u)  \le a_{\delta,L}(u)    \le (1+\eta) a_{\delta,K}(u), $$
whence $$(1-\eta)K_{\delta}\subseteq  L_{\delta}\subseteq (1+\eta)K_{\delta}.$$
Since $rB_2^n\subseteq K_{\delta}\subseteq RB_2^n$, it then follows  that, given $\rho>0$, for $\eta>0$ small enough, one has 
$d_H(L_{\delta},K_{\delta})\le \rho.  $
$\square$
\vskip 3mm
As a corollary, we obtain new affine invariant points.
\vskip 2mm
\noindent
\begin{corollary}\label{continuity}
Let $0<\delta<\big(\frac{n}{n+1}\big)^n$ and let $p: {\cal K}_n\to \mathbb{R}^n$  be an affine invariant point.  Then $K\to
p(K_{\delta})$ is also an affine invariant point. In particular, for the centroid $g$,  $K\mapsto g(K\setminus K_{\delta})$ is an affine
invariant point. \end{corollary}
\vskip 2mm 
\noindent
{\bf Proof.} Affine invariance follows  from  Remark (i) after Definition \ref{ais} and continuity from Proposition \ref{cont:Kdelta}. The second statement  
follows now from the trivial identity $$g(K) =\frac{|K_\delta |}{|K|} g(K_{\delta})+  \frac{| K \setminus K_\delta |}{|K|} g(K\setminus K_{\delta}), $$
which gives $$g(K\setminus K_{\delta})= \frac{|K|}{| K \setminus K_\delta |} \  g(K)   -  \frac{|K_\delta|}{| K \setminus K_\delta |} \  g(K_{\delta}),$$
as an affine combination of continuous affine invariant points. $\square$
\vskip 3mm
The next lemma is key for many of the proofs that will follow.
\vskip 2mm
\noindent
\begin{lemma} \label{gincap}
Let $m \geq n+1$ and for $1 \leq i \leq m$,  let $v_i \in \mathbb{R}^n$ be the vertices of 
a polytope $P$ in $\mathcal{K}_n$. For all $\varepsilon >0$ there exists  $z \in P$ with 
$\|v_{1}-z\|\leq\varepsilon$ and $0 < r  \leq \varepsilon$ such that $B^n_2(z,r) \subset P$
and if 
$
K = \conv\left(B^n_2(z,r), v_2, \dots , v_m\right) 
$, 
then $K$ satisfies 
\par
(i) $ K \subseteq P$, $v_2, \dots , v_m$ are extreme points of $K$  and $d_H(K,P) \leq \varepsilon$.
\par
(ii) For sufficiently small $\delta$, $
\|v_1- g\left(K \setminus K_{\delta } \right)\| \leq2 \varepsilon  $
\end{lemma}
\vskip 2mm
\noindent
{\bf Proof.}
There exists 
a hyperplane $H$ that  striclty separates $v_1$ and $\{v_2, \dots, v_m\}$,  such that
 for all $x \in H^- \cap P$ we have that $\|x-v_1\| < \varepsilon$.
Let $z \in \inte(H^-) \cap \inte(P)$. Then  there exists $0 < r  \leq \varepsilon $ such that $B^n_2(z,r) \subseteq H^-\cap P$. Let 
$K = \conv\left(B^n_2(z,r), v_2, \dots , v_m\right) .$
\par
\noindent
(i)  By construction of $K$, 
$v_2 \dots, v_m$ are extreme points of $K$. 
 Also, for all   $x \in K \cap H^-$,  one has $\|x-v_1\| < \varepsilon$. Therefore $K \subseteq  P\subseteq K + \varepsilon B^n_2$ and thus $d_H(K,P) \leq \varepsilon$.
\par
\noindent
(ii) We have for
$\delta>0$, 
\begin{eqnarray}\label{split}
g\left(K \setminus K_{\delta } \right) &=& \frac{\left| \left( K \cap H^+\right) \setminus K_{\delta }\right| }{\left|K \setminus K_{\delta } \right|}  g \left( \left( K \cap H^+\right) \setminus K_{\delta  }\right)\nonumber \\
&+& \frac{\left| \left( K \cap H^-\right) \setminus K_{\delta  }\right| }{\left|K \setminus K_{\delta} \right|}  g \left( \left( K \cap H^-\right) \setminus K_{\delta }\right).
\end{eqnarray}
Since $g \left( \left( K \cap H^-\right) \setminus K_{\delta  }\right) \in \inte \left( K\right)  \cap \inte\left( H^-\right)$, one has 
\begin{equation*}
\|v_1 -  g \left( \left( K \cap H^-\right) \setminus K_{\delta  }\right)\| \leq \varepsilon.
\end{equation*}
Observe that $\partial K$ contains a cap of 
$\partial B(z,r)$, so that 
$$C= \int _{\partial K}  \kappa_K^{\frac{1}{n+1}} d\mu_K >0.$$
By Theorem \ref{SW1990},  one has  for $\delta$ sufficiently small,
\begin{equation}\label{est22}
|K \setminus K_{\delta  }| =|K|-|K_{\delta}| \geq  
\frac{C}{2c_n} \left(\delta |K|\right)^{\frac{2}{n+1}}.
\end{equation}
Let $R=\max\{\|x\|: x\in P\}$.
As the Gauss curvature is equal to $0$ everywhere on the boundary $\partial \left(K \cap H^+\right)$, again by  Theorem \ref{SW1990},   one has for 
sufficiently small  $\delta$, 
$$
c_n \frac{\left|K \cap H^+\right| - \left|\left(K \cap H^+\right)_{\frac{\delta |K|}{|K\cap
H^+|} }\right|}{\left(\delta |K|\right)^\frac{2}{n+1}} \leq \frac{C\varepsilon}{4R }. $$
As 
$$
\left( K \cap H^+\right) \setminus K_{\delta } 
\subseteq \left( K \cap H^+\right) \setminus \left(K  \cap H^+\right)_{\frac{\delta |K|}{|K\cap H_2^+|}},
$$
we get 
\begin{equation} \label{est1}
\left | \left(K \cap H^+\right) \setminus K_{\delta }\right| \leq \frac{C\varepsilon}{4Rc_n} \       \left(\delta |K|\right)^\frac{2}{n+1}.
\end{equation}
It folllows from  (\ref{est22}) and (\ref{est1}) that for $\delta $ small enough one has
\begin{equation} \label{est3}
\frac{\left| \left( K \cap H^+\right) \setminus K_{\delta }\right| }{\left|K \setminus K_{\delta } \right|} \leq \frac{\varepsilon}{R} .
\end{equation} 
We get thus from (\ref{split}) and (\ref{est3})
\begin{eqnarray*} 
&& \| g\left(K \setminus K_{\delta } \right) - g \left( \left( K \cap H^-\right) \setminus K_{\delta  }\right)\|   \\
&&= \frac{\left| \left( K \cap H^+\right) \setminus K_{\delta }\right| }{\left|K \setminus K_{\delta } \right|}
\| g \left( \left( K \cap H^+\right) \setminus K_{\delta  }\right) 
 - g \left( \left( K \cap H^-\right) \setminus K_{\delta  }\right)\| \\
&&\leq \frac{\left| \left( K \cap H^+\right) \setminus K_{\delta }\right| }{\left|K \setminus K_{\delta } \right|} \left(\| g \left( \left( K \cap H^-\right) \setminus K_{\delta  }\right) \| + \| g \left( \left( K \cap H^+\right) \setminus K_{\delta  }\right) \| \right) \\
 &&\leq \frac{\varepsilon}{2R} (R+R)\leq  \varepsilon
\end{eqnarray*}
Altogether, 
$$\|v_1 -  g\left(K \setminus K_{\delta } \right)\| \leq  \| v_1 - g \left( \left( K \cap H^-\right) \setminus K_{\delta  }\right)\| + \| g\left(K \setminus K_{\delta } \right) - g \left( \left( K \cap H^-\right) \setminus K_{\delta  }\right)\| \leq 2 \varepsilon.
$$
 $\square$
\vskip 3mm
\noindent
\vskip 3mm
\subsection{Proof of Theorem \ref{unendlich}: $\mathfrak{P}_n$ is infinite dimensional.}
Here,  we answer in the negative Gr\"unbaum's question whether there exists a finite basis for $\mathfrak{P}_n$, 
 i.e. affine invariant points 
$p_i\in \mathfrak{P}_n$, $1 \leq i \leq l$,  such that every  $p\in \mathfrak{P}_n$ can be written as $$p=\sum_{i=1}^l \alpha_i p_i, \hskip 3mm\text{ with } \alpha_i \in \mathbb{R}\hskip 2mm\text{ and  }\sum_{i=1}^l \alpha_i =1.$$
\par 
Recall that $\mathfrak{P}_n$ is an affine subspace of $C(\mathcal{K}_n, \mathbb{R}^n)$, the continuous functions on $\mathcal{K}_n$ with values in $  \mathbb{R}^n$ and that we  denote by $V\mathfrak{P}_n$ the subspace parallel to $\mathfrak{P}_n$. Thus, with  the centroid $g$, 
\begin{equation}\label{vektorraum}
V\mathfrak{P}_n = \mathfrak{P}_n - g.
\end{equation}
The dimension of $\mathfrak{P}_n$ is the dimension of $V\mathfrak{P}_n$. We introduce a norm on $V\mathfrak{P}_n$,  
\begin{equation}\label{vektorraumnorm}
\|v\|=\sup_{K\in\mathcal{K}_{n},
B_{2}^{n}\subseteq K\subseteq n B_{2}^{n}}\|v(K)\|, \hskip 8mm \text{for} \  v \in V\mathfrak{P}_n.
\end{equation}
Observe  that the set $
\{K\in\mathcal K_{n}: B_{2}^{n}\subseteq K \subseteq n B_{2}^{n}\}
$
is a compact subset of $(\mathcal K_{n}, d_H)$. Therefore
(\ref{vektorraumnorm}) is well defined and it is a norm: $v=p-g \neq 0$ implies 
 that there is $C$ with $v(C)\ne0$. By John's theorem (e.g., \cite{NicoleBuch}),  there is 
an affine, invertible map $T$ with
$B_{2}^{n}\subseteq T(C)\subseteq nB_{2}^{n}$. Thus,
\begin{eqnarray*}
v(T(C))
=(p-g)(T(C))
=p(T(C))
-g(T(C))
=T(p(C))
-T(g(C)).
\end{eqnarray*}
Since $T=S+x_{0}$, where $S$ is a linear map,
$$
v(T(C))
=S(p(C)
-g(C))
=S((p-g)(C))\ne0.
$$
Hence 
$$
\|v\|\geq\|v(T(C))\|>0.
$$

\vskip 3mm
For the proof of Theorem \ref{unendlich} and Theorem \ref{dicht}, we will make use of the following
theorem by Sch\"utt   and Werner \cite{SchuettWerner1990}.  There, $\mu_K$ is the usual surface measure on $\partial K$ and for $x\in\partial K$, $\kappa(x)$ is the generalized Gauss
curvature at $x$, which is defined $\mu_K$ almost everywhere.
\begin{theorem} \cite{SchuettWerner1990} \label{SW1990}
Let $K$ be a convex body in $\mathbb{R}^n$. Then, if  $c_n= 2  \left(\frac{|B^{n-1}|}{n+1}\right)^{\frac{2}{n+1}}$, one has
$$
c_n \lim_{\delta \rightarrow 0} \frac{|K| - |K_\delta|}{\left(\delta |K|\right)^\frac{2}{n+1}} = \int_{\partial K} \kappa^\frac{1}{n+1}(x)\  d\mu_K(x).
$$
\end{theorem}
\vskip 3mm
\noindent
{\bf Proof of Theorem \ref{unendlich}.}
We show that the closed unit ball of $V\mathfrak{P}_{n}$ is not compact.
For $K\in {\mathcal K}_n$ and $\delta>0$,  let $K_{\delta}$ be the convex floating
body of $K$. Let $g$ be the centroid and let $g_{\delta}:\mathcal K_{n}\rightarrow\mathbb R^{n}$ be the affine invariant point  given by
$$
g_{\delta}(K)=g(K\setminus K_{\delta}).
$$
The set of vectors \{$v_{\delta}=g_{\delta}-g$: $\delta>0\}$ is bounded. Indeed,
since $g(K)\in K$ and $g_{\delta}(K)\in K$
$$
\|v_{\delta}\|
\leq
\sup_{K\in\mathcal{K}_{n},
B_{2}^{n}\subseteq K\subseteq n B_{2}^{n}}(\|g(K)\|+\|g_{\delta}(K)\|)
\leq2n.
$$
The sequence $v_{\frac{1}{j}}$, $j\geq1$, does not have a convergent
subsequence:  We show that for all $N$ there are $\ell\geq m\geq N$ and $K\in\mathcal K_{n}$ with $B_{2}^{n}\subseteq K\subseteq n B_{2}^{n}$ such that
$$
\left\|v_{\frac{1}{\ell}}(K)-v_{\frac{1}{m}}(K)\right\| \geq \frac{1}{10}.
$$
As $K$ we choose the union of the cylinder $D=[-1,1]\times B_2^{n-1}$ and a cap of a Euclidean ball, 
\begin{equation}\label{kh}
K(h)=D \cup\ C(h),
\end{equation}
where, with $e_1=(1,0,\dots,0)\in \R^n$,
$$
C(h)=\left( \frac{h^{2}+2h-1}{2h}e_{1}+ \frac{1+h^{2}}{2h}B_2^n\right)
\cap \{x=(x_1,\dots,x_n)\in \R^n: x_1\ge 1\}.
$$
As $h\to 0$, $K(h)\to D$ and, by Corollary \ref{continuity}, 
$
g_{\frac{1}{m}}(K(h))\to
g_{\frac{1}{m}}(D)=0.
$
Thus there exists 
$h_{0}>0$ such that  
$$
\|g_{\frac{1}{m}}(K(h))\|
=\left\|g\left(K(h)\setminus(K(h))_{\frac{1}{m}}\right)\right\|
\leq\frac{1}{10}, \hskip 2mm \hbox { for all $h\le h_0$}.
$$
Now we show that we can choose $\ell$ sufficiently big so that
\begin{equation}\label{infinite1}
\|g_{\frac{1}{\ell}}(K(h))\|\geq\frac{1}{2}.
\end{equation}
We apply the same argument as in the proof of Lemma \ref{gincap}.
Let $H$ be the hyperplane such that
$$
K(h)\cap H^{-}=C(h)
\hskip 10mm\mbox{and}\hskip 10mm
K(h)\cap H^{+}=D.
$$
Then as in (\ref{split}),
$$
g_{\frac{1}{\ell}}(K(h))
= \frac{\left| D
 \setminus K(h)_{\frac{1}{\ell} }\right| }{\left|K(h) \setminus K(h)_{\frac{1}{\ell}} \right|}  
 g \left( D \setminus K(h)_{\frac{1}{\ell}  }\right)\nonumber \\
+\frac{\left| C(h) \setminus K(h)_{\frac{1}{\ell}  }\right| }{\left|K(h) \setminus K(h)_{\frac{1}{\ell}} \right|}  
g \left(  C(h)\setminus K(h)_{\frac{1}{\ell} }\right).
$$
Since 
$$
g(D\setminus K(h)_{\frac{1}{\ell}})\in D
\hskip 10mm
\mbox{and}
\hskip 10mm
g(C(h)\setminus K(h)_{\frac{1}{\ell}})\in C(h)
$$
we get
$$
\|g(D\setminus K(h)_{\frac{1}{\ell}})\|
\leq \sqrt{2}
\hskip 10mm\mbox{and}\hskip 10mm
\|g(C(h)\setminus K(h)_{\frac{1}{\ell}}\|\geq1.
$$
Therefore, by triangle inequality,
\begin{equation}\label{gincap2}
\|g_{\frac{1}{\ell}}(K(h))\|
\geq\frac{\left| C(h) \setminus K(h)_{\frac{1}{\ell}  }\right| }{\left|K(h) \setminus K(h)_{\frac{1}{\ell}} \right|} 
-\sqrt{2}\frac{\left| D
 \setminus K(h)_{\frac{1}{\ell} }\right| }{\left|K(h) \setminus K(h)_{\frac{1}{\ell}} \right|}  .
\end{equation}
By Theorem \ref{SW1990},  we get as in (\ref{est22}), for $\ell$ large enough, with $ \alpha(h)=\int _{\partial K(h)}  \kappa_{K(h)}^{\frac{1}{n+1}} d\mu_{K(h)} $, 
\begin{equation*}
\left|K(h) \setminus K(h)_{\frac{1}{\ell} }\right| 
\geq \frac{\left(\frac{1}{\ell} |K(h)|\right)^\frac{2}{n+1}}{2 c_n} \int _{\partial K(h)}  \kappa_{K(h)}^{\frac{1}{n+1}} d\mu_{K(h)} 
= \frac{\left(\frac{1}{\ell} |K(h)|\right)^\frac{2}{n+1}}{2 c_n}  \alpha(h).
\end{equation*}
Also by Theorem \ref{SW1990},  we get as in (\ref{est1}),
\begin{equation*} \label{}
\left|D\setminus K(h)_{\frac{1}{\ell}}\right|
=
\left | \left(K(h) \cap H^+\right) \setminus K(h)_{\frac{1}{\ell} }\right|
 \leq \frac{\varepsilon}{c_n} \   \alpha(h) \    \left(\frac{1}{\ell} |K(h)|\right)^\frac{2}{n+1}.
\end{equation*}
Now we finish the proof as in Lemma \ref{gincap}. $\Box$

\subsection {Proof of Theorem \ref{Pn=Fn}.  } 

It was also asked by Gr\"unbaum \cite{Gruenbaum1963} if 
for every $K \in \mathcal{K}_n$, 
$$\mathfrak{P}_n(K) = \mathfrak{F}_n(K),$$
where
$\mathfrak{F}_n(K)= \{ x \in \mathbb{R}^n: Tx=x, \text{  for all affine  } T 
\text{ with  } TK=K\}$. Observe that it is clear that $ \mathfrak{P}_n(K) \subseteq \mathfrak{F}_n(K)$.
We will prove  that $\mathfrak{P}_n(K) = \mathfrak{F}_n(K)$, 
if $ \mathfrak{P}_n(K) $ is $(n-1)$-dimensional. To do so, we, again, 
first  need  to define  new affine invariant set mappings.
\par
Actually, in the proof of  Theorem \ref{Pn=Fn} we show that the group of isometries of
$K$ equals
$$
\{I_n,S\}=\{T:\R^n\to \R^n \hbox{ affine one to one},  \hskip 2mm TK=K\},
$$
where $S$ is  reflection about  a hyperplane, i.e. $S: \mathbb{R}^n \rightarrow \mathbb{R}^n$ is bijective and there is a hyperplane $H$
and a direction $\xi \notin H$ such that $S(h + t \xi)= h-t \xi$ for all $h \in H$.
\vskip 3mm

\begin{lemma}
(i)
Let $p \in {\mathfrak P}_n$  and let $g$ be the centroid. For $0 < \varepsilon  < 1$, define 
$$
A_{p,\varepsilon}(K)
=\left\{x\in K\left|\langle x,p((K-g(K))^{\circ})\rangle 
\geq (1 - \varepsilon)     \sup_{y\in K}\langle y,p((K-g(K))^{\circ}) \rangle \right.\right\}.
$$
Then $A_{p,\varepsilon}:\mathcal K_{n}\rightarrow\mathcal K_{n}$ is an affine invariant set map.
\vskip 2mm
\noindent
(ii) Let $p\in {\mathfrak P}_n$  and let $q\in {\mathfrak P}_n$ be  proper.
Then $A_{q,p,\varepsilon}:\mathcal K_{n}\rightarrow\mathcal K_{n}$ given by
$$
A_{q,p,\varepsilon}(K)
=\left\{x\in K\left|\langle x, p((K-g(K))^{\circ})\rangle
\geq   \left(1- \varepsilon\right)  \langle q(K),p((K-g(K))^{\circ})\rangle \right.\right\}
$$
is an affine invariant set map.
\end{lemma}
\vskip 3mm

Observe that
$0\in{\mathfrak P}_n((K-g(K))^{\circ})$ since $0$ is the Santal\'o point of 
$(K-g(K))^{\circ}$. Therefore, ${\mathfrak P}_n((K-g(K))^{\circ})$ is a subspace
of $\mathbb R^{n}$.
\vskip 3mm

\noindent
{\bf Proof.} Let $T$ be an invertible, affine map and $T=S+a$ its
decomposition in a linear map $S$ and a translation $a$.
Then for any  convex body $C$ that contains $0$ in its interior, 
        $$
(S( C))^{\circ}=S^{*-1}(C^{\circ}).
$$
Moreover,
\begin{eqnarray*}
p((T(K)-g(T(K)))^{\circ})
&=&p((S(K-g(K)))^{\circ})  \\
&=&p(S^{*-1}((K-g(K))^{\circ}))
=S^{*-1}(p((K-g(K))^{\circ})).
\end{eqnarray*}
Since $S^{*-1*}=S^{-1}$, 
\begin{eqnarray*}
&&\hskip -4mm A_{p,\varepsilon}(T(K))   = \\
&&\hskip -4mm=\left\{x\in T(K)\left|\langle x,p((T(K)-g(T(K)))^{\circ}) \rangle
\geq (1- \varepsilon) \sup_{y\in T(K)}\langle y,p((T(K)-g(T(K)))^{\circ})\rangle \right.\right\}  \\
&&\hskip -4mm=\left\{x\in T(K)\left|\langle S^{-1}x,p((K-g(K))^{\circ})\rangle
\geq (1- \varepsilon) \sup_{y\in T(K)} \langle S^{-1}y,p((K-g(K))^{\circ}) \rangle\right.\right\}
\end{eqnarray*}
and one verifies easily that $A_{p,\varepsilon}(T(K)) =T(A_{p,\varepsilon}(K))$.
Please note that $A_{p,\varepsilon}(K)$ is convex, compact and nonempty.
$\Box$
\vskip 3mm

\begin{lemma}\label{dimdim}
Let $K \in \mathcal{K}_n$  and let 
$P:\mathbb R^{n}\rightarrow\mathbb R^{n}$ be the orthogonal
projection onto $\mathfrak P_{n}((K-g(K))^{\circ})$. Then the
restriction of $P$ to the subspace $\mathfrak P_{n}(K-g(K))$
is an isomorphismn between $\mathfrak P_{n}(K-g(K))$ and
$\mathfrak P_{n}((K-g(K))^{\circ})$.
\par
In particular,
$$
\operatorname{dim}(\mathfrak P_{n}(K-g(K)))
=\operatorname{dim}(\mathfrak P_{n}((K-g(K)))^{\circ})
$$
\end{lemma}
\vskip 3mm

\noindent
{\bf Proof.}
On the hyperplane
$\mathfrak P_{n}((K-g(K))^{\circ})$, $P(K-g(K))$ has an interior point. This holds because otherwise,
by Fubini, $\operatorname{vol}_{n}(K)=0$.
\par
Let 
$k=\operatorname{dim}(\mathfrak P_{n}((K-g(K)))^{\circ})$.
We choose $u_{1}\in\mathfrak P_{n}((K-g(K))^{\circ})$. Then 
$g(A_{u_{1},\varepsilon_{1}})$ is a proper affine invariant point. Now we choose
$u_{2}\in \mathfrak P_{n}((K-g(K))^{\circ})$ that is orthogonal to
$P(g(A_{u_{1},\varepsilon_{1}}))$. Then 
$P(g(A_{u_{1},\varepsilon_{1}}))$ and $P(g(A_{u_{2},\varepsilon_{2}}))$ are linearly
independent.
\par
Eventually,
$$
P(g(A_{u_{1},\varepsilon_{1}})),\dots,P(g(A_{u_{k},\varepsilon_{k}}))
$$
are linearly independent, and therefore
$$
g(A_{u_{1},\varepsilon_{1}}),\dots,g(A_{u_{k},\varepsilon_{k}})
$$
are linearly independent. Therefore,
$$
\operatorname{dim}(\mathfrak P_{n}((K-g(K)))^{\circ})
\leq\operatorname{dim}(\mathfrak P_{n}(K-g(K))).
$$
Now we interchange the roles of $\mathfrak P_{n}(K-g(K))$ and
$\mathfrak P_{n}((K-g(K))^{\circ})$ and get the inverse inequality.
\par
Let $Q$ denote the restriction of $P$ to the subspace
$\mathfrak P_{n}(K-g(K))$.  
$g(A_{u_{1},\varepsilon_{1}}),\dots,g(A_{u_{k},\varepsilon_{k}})$
is a basis of $\mathfrak P_{n}(K-g(K))$ and 
$P(g(A_{u_{1},\varepsilon_{1}})),\dots,P(g(A_{u_{k},\varepsilon_{k}}))$
is a basis of $\mathfrak P_{n}((K-g(K)))^{\circ}$. $Q$ is a bijection
between the two bases, thus $Q$ is an isomorphism.
$\Box$
\vskip 3mm

\begin{lemma}\label{SurProp}
Let $K\in\mathcal K_{n}$. Then for every point $x$ from the relative interior
of $K\cap\mathfrak P_{n}(K)$ there is a proper affine invariant point $q$ with
$q(K)=x$.
\end{lemma}
\vskip 2mm
\noindent
{\bf Proof.} We use the same notation as in Lemma \ref{dimdim} and its proof.
We may assume that $g(K)=0$. Suppose that there is an interior point
$x$ of $\mathfrak P_{n}(K)\cap K$ in the hyperplane $\mathfrak P_{n}(K)$
for which there is no proper affine invariant point $q$ with $q(K)=x$.
The set 
$$
\{p(K)| p\hskip 1mm \mbox{is a proper affine invariant point}\}
$$
is convex. $P:\mathbb R^{n}\rightarrow\mathbb R^{n}$ is the orthogonal
projection onto $\mathfrak P_{n}(K^{\circ})$. Then 
$P(\mathfrak P_{n}(K)\cap K)$ is a convex set in the hyperplane 
$\mathfrak P_{n}(K^{\circ})$. Since $P$ is an isomorphism between the
hyperplanes $\mathfrak P_{n}(K^{})$ and $\mathfrak P_{n}(K^{\circ})$ we have
$$
P(x)\notin P(\{p(K)| p\hskip 1mm \mbox{is a proper affine invariant point}\}).
$$ 
Moreover, $P(x)$ is an interior point of $P(\mathfrak P_{n}(K)\cap K)$.
By the Hahn-Banach theorem there is $u\in\mathfrak P_{n}(K^{\circ})$
such that for all proper affine invariant points $p$ we have
$$
\langle u,x\rangle \geq \langle u,P(p(K)) \rangle.
$$
On the other hand, there is an affine invariant point $q$ with $q(K^{\circ})=u$. Then
$g\circ A_{u,\langle u,x \rangle}$ is a proper affine invariant point
with
$$
\langle u,x \rangle < \langle u,g\circ A_{q,\langle u,x \rangle } \rangle, 
$$
which is a contradiction.
$\Box$
\vskip 3mm

\begin{lemma} \label{Lemma-Pn=Fn}
Let $K \in \mathcal {K}_n$ and suppose that
$\operatorname{dim}(\mathfrak P_{n}(K))=n-1$. 
Then $S:\mathbb R^{n}\rightarrow\mathbb R^{n}$ with
$$
S(y+x)=y-x
$$
for all $y\in \mathfrak P_{n}(K-g(K))$ and 
$x\in \mathfrak P_{n}((K-g(K))^{\circ})^{\perp}$ is a linear map such that 
$$
S(K-g(K))=K-g(K).
$$
\end{lemma}
\vskip 2mm
 \noindent
{\bf Proof.}
By Lemma \ref{dimdim}, the orthogonal projection onto
$\mathfrak P_{n}((K-g(K))^{\circ})$ restricted to $\mathfrak P_{n}(K-g(K))$
is an isomorphism. Therefore, 
$$
\mathbb R^{n}
=\mathfrak P_{n}(K-g(K))\oplus\mathfrak P_{n}((K-g(K))^{\circ})^{\perp}.
$$
By Lemma \ref{SurProp}
for every $y\in\mathfrak P_{n}(K-g(K))\cap\operatorname{int}(K)$ there is a
proper affine invariant point $q$ with $y=q(K)$. Let $u_{1},\dots,u_{n-1}$
be an orthonormal basis in $\mathfrak P_{n}((K-g(K))^{\circ})$.
The map $A_{\varepsilon}:\mathcal K_{n}\rightarrow\mathcal K_{n}$ defined by
$$
A_{\varepsilon}(K)
=\bigcap_{i=1}^{n-1}\{x\in K| \langle q(K),u_{i}\rangle -\varepsilon
\leq\langle x,u_{i} \rangle \leq\langle q(K),u_{i} \rangle +\varepsilon\}
$$
is an affine invariant set map. As $ q(K)$ is an interior point of 
$A_{\varepsilon}(K)$, $A_{\varepsilon}(K) \in \mathcal{K}_n$. Moreover,
$$
\lim_{\varepsilon\to0}A_{\varepsilon}(K)
=K\cap(q(K)+\mathfrak P_{n}((K-g(K))^{\circ})^{\perp})
$$
in the Hausdorff metric.
 $g\circ A_{\varepsilon}$ is a proper affine invariant
point. Since all affine invariant points are elements of $\mathfrak P_{n}(K)$ 
$$
\lim_{\varepsilon\to0}(g\circ A_{\varepsilon})(K)
=q(K).
$$
On the other hand, $q(K)$ is the midpoint of
$K\cap(q(K)+\mathfrak P_{n}((K-g(K))^{\circ})^{\perp})$.
$\Box$
\vskip 3mm
\noindent
{\bf Proof of Theorem \ref{Pn=Fn}.} 
Theorem \ref{Pn=Fn} now follows immediately from Lemma \ref{Lemma-Pn=Fn}. 
Indeed, Lemma \ref{Lemma-Pn=Fn} provides a map $T= S - S(g(K)) + g(K)$  with $T(K)=K$
and such that for all $z\in \mathfrak {P}_{n}(K)$ and for all 
$x\in \mathfrak {P}_{n}((K-g(K))^{\circ})^{\perp}$, 
$$T(z+x) = z-x.
$$
Consequently, if $w\notin \mathfrak {P}_{n}(K)$, then $T(w)\neq w$, which means that the complement of $\mathfrak {P}_{n}(K)$
is contained in the complement of $\mathfrak {F}_{n}(K)$.
\vskip 3mm

\noindent
{\bf Remark.} As a byproduct of the preceding results, it can be proved that  if $K\in \K$
satisfies  ${\mathfrak P}_n(K)=\R^n$ and if
${\mathfrak S}_n(K) =\{A(K): A\in {\mathfrak S}_n\}$, 
then ${\mathfrak S}_n(K)$ is dense in $\K$. 
It might be conjectured that for general $K\in \K$,  
${\mathfrak S}_n(K)$ is dense in $\{C\in \K: \ {\mathfrak F}_n(C)\subseteq
{\mathfrak F}_n(K)\}.$

\subsection{Proof of Theorem  \ref{dicht}.}
In this subsection we show that the set of all $K$ such that $\mathfrak{P}_n(K)=\mathbb{R}^n$,  is dense in $\mathcal{K}_n$ and consequently
the set of all $K$ such that $\mathfrak{P}_n(K)=\mathfrak{F}_n(K)$ is dense in $\mathcal{K}_n$.
A further corollary is that,  for  every  $k \in \mathbb{N}$, $0 \leq k \leq n$, there exists a convex body $Q_k$
such that $\mathfrak{P}(Q_k)$ is a $k$-dimensional affine subspace of $\mathbb{R}^n$.
\par
It is relatively easy to construct   examples of convex bodies $K$ in the plane such that  $\mathfrak{P}_n(K)=\mathbb{R}^2$.
To do so in higher dimensions is more involved and we present  a construction in  the proof of Theorem \ref{dicht} below.
First, we will briefly mention two examples in the plane.
\vskip 3mm\noindent
{{\bf Example 1.}
Let $S$ be a regular simplex in the plane and let ${\cal J}(S)$ be the ellipsoid of maximal  area inscribed in $S$.
We show in the section  below that the center $j(S)$ of ${\cal J}(S)$  is an affine invariant point. We can assume that ${\cal J}(S)= B^2_2$, the Euclidean ball centered at 
$0$ with radius $1$. Then e.g. $S = \conv \left( (-1, -\sqrt{3}), (-1, \sqrt{3}), (2,0)\right)$.
\par
Let  $0 <\lambda <1$ be such that $H((1+ \lambda) e_1, e_1) \cap \mbox{int}(S) \neq \emptyset$ and  consider the convex body $S_1=S \cap H^+((1+ \lambda) e_1, e_1)$ obtained from $S$ by cutting of a cap from $S$. Then still $j(S_1)=0$ but the center 
of gravity has moved to the left of $0$.  Next, 
let  $\gamma >0$ be such that $H((1+ \gamma) u,u) \cap \mbox{int}(S_1) \neq \emptyset$, where $u=\frac{ (-1, \sqrt{3})}{2}$ and  consider the convex body $S_2=S_1 \cap H^+((1+ \gamma u), u)$ obtained from $S_1$ by cutting of a cap from $S_1$. Then still $j(S_2)=0$ but the center 
of gravity $g(S_2)$ of $S_2$ has moved and it is different from the Santal\'o point $s(S_2)$ of $S_2$. $j(S_2)$, 
$g(S_2)$ and $s(S_2)$ are three affinely independent points of $\mathbb{R}^2$, hence span $\mathbb{R}^2$.
\vskip 2mm
\noindent
{{\bf Example 2.}
Let $S$ be the  equilateral triangle in the plane centered at $0$ of Example 1 with vertices 
$a=(2,0)$,  $b= (-1, \sqrt{3})$ and  $c= (-1, -\sqrt{3})$.  Then, as noted in Example 1,  $B^2_2$ is the John ellipse ${\cal J}(S)$ of 
$S$. Let $b_1 , c_1$ be two points on the segments
$[a,b]$ and $[a,c]$, such that the segment $[b_1,c_1]$ does not intersect $B^2_2$. Then $B^2_2$ is still the John ellipse of the quadrangle $\conv \left(b, b_1, c_1, c \right)$. 
 Now the L\"owner ellipse
${\cal L}(T)$ of
the triangle $T=\conv \left(b_1, b, c\right)$ is centered at $\frac{1}{3}(b_1+b+c)\not=0$, if $b_1\not=a$. ${\cal L}(T)$  intersects  the segment $[a,c]$ at $c$ and at some point $c'$. When $b_1\to a$, one has ${\cal L}(T) \to 2 B^2_2 $ and thus $c'\to a$. So we may
choose $b_1$ such that $[b_1,c']$ does not meet $B^2_2$, and  thus  for some $c''\in [a,c]$, $[b_1, c_1]$ does not meet
$B^2_2$ for any $c_1\in [c',c'']$. Finally,  let $P(c_1)$ be the quadrangle   $P(c_1)=\conv \left(b, b_1, c_1, c\right)$, with $c_1\in [c',c'']$. Since $b_1,  b, c$ and $c_1$ are the vertices of $P(c_1)$ and $c_1\in {\cal L}(T)$, ${\cal L}(T)$ is also  the L\"owner ellipsoid ${\cal L}(P(c_1))$of $P(c_1)$. Altogether, 
\par
The John ellipse of $P(c_1)$ is $B^2_2 $ which is centered at  $0$, so that the affine invariant point $j(P(c_1))=0$.
\par
The L\"owner ellipse of $P(c_1)$ is centered at $\frac{1}{3}(b_1+b+c)$, so that the affine invariant point $l((P(c_1))= \frac{1}{3}(b_1+b+c) \neq 0$.
\par
An easy computation shows that the centroid of $P(c_1)$ moves on an hyperbola when $c_1$ varies in $[c',c'']$.
\par
So, in general,  these three points are not on line.  $\square$
\vskip 2mm
\noindent
{\bf Proof of Theorem \ref{dicht}.}
The set of $n$-dimensional polytopes  is dense in $\left(\mathcal{K}_n, d_H\right)$. Let $P$ be a polytope and let $\eta >0$ be given. Then it is enough to show that 
there exists a convex body $Q$ with $d_H(P,Q) < \eta$ and such that
$\mathfrak{P}_n(Q) = \mathbb{R}^n$. 
\par
We  describe the idea of the proof. For a properly constructed convex body $Q$ we will construct  $\Delta_i\in{\mathfrak P}_n$, $1 \leq i \leq n+1$,  in  such a way that the 
$\Delta_i(Q)$ are affinely independent.
\par
The construction of such a $Q$ is  done  inductively: we first construct  $Q_1$ very near $P$ and such that $\Delta_1(Q_1)$ is near an extreme point $v_1$  of $Q_1$. Then we construct  $Q_2$ very near $Q_1$ and $P$ and such that $\Delta_1(Q_2)$ is near the  extreme point $v_1$ of $Q_2$
and $\Delta_2(Q_2)$ is near an  extreme point  $v_2 \neq v_1$ of $Q_2$.
\par
Let $P = \conv\left(v_1, \dots, v_m\right)$ be a polytope with non-empty interior and with $m$ vertices, $m \geq n+1$.  
We  pick $n+1$  affinely independent vertices of $P$. We can assume that these are $v_1, \dots, v_{n+1}$. 
Let $0 < \eta_1 < \frac{\eta}{n+2}$ be given. By Lemma \ref{gincap}, there exists $z_1 \in P$, $\|v_{1}-z_{1}\|\leq\eta_1$,  and $0 < r_1  \leq \eta_1$ such that $B^n_2(z_1,r_1) \subseteq P$
and such that 
$$
Q_1 = \conv\left(B^n_2(z_1,r_1), v_2, \dots , v_m\right) 
$$
has  $v_2, \dots , v_m$ as  extreme points,    
\begin{equation}\label{dq1p}
d_H(Q_1,P) \leq \eta_1,
\end{equation}
and for sufficiently small $\delta_1$, 
\begin{equation}\label{v1-q1}
\|v_1- g\left(Q_1 \setminus (Q_1)_{\delta_1 } \right)\| \leq 2\eta_1.
\end{equation}
We let $\varepsilon_1< \eta_1$ and choose an $\varepsilon_1$-net $\mathcal{P}_{\varepsilon_1}$ on $\partial \left(B^n_2(z_1,r_1) \right)$ and put
$$
P_1=\conv\left( \mathcal{P}_{\varepsilon_1}, v_2, v_3, \dots, v_m\right).
$$
Then $P_1 \subseteq Q_1 \subseteq P$ and  $d_H(P_1,Q_1) \leq \varepsilon_1 < \eta_1$.
By Corollary \ref{continuity}, for  a given $K\in \K$, for a given $0<\delta<\big(\frac{n}{n+1}\big)^n$ and $\varepsilon>0$,  there exists $\gamma(K, \delta, \varepsilon)$
such that if 
\begin{equation}\label{cor1}
d_H(K,L)< \gamma(K, \delta, \varepsilon), \ \text{  for $L\in \K$, } \   \text{ then } \  
\left\| g(K\setminus K_{\delta}) - g(L\setminus L_{\delta})\right\| <  \varepsilon.
\end{equation}
As $d_H(P_1,Q_1) \leq \varepsilon_1$, we get that
$$
\left\| g(Q_1\setminus (Q_1)_{\delta_1}) - g(P_1\setminus (P_1)_{\delta_1})\right\| <  \eta_1,
$$
if we choose  in addition $\varepsilon_1$ such that $\varepsilon_1 < \gamma(Q_1, \delta_1, \eta_1)$.
Thus, together with (\ref{v1-q1}), 
\begin{equation}\label{v1p1}
\left\| v_1- g(P_1\setminus (P_1)_{\delta_1})\right\| \leq 3 \  \eta_1.
\end{equation}
Observe that $v_{2},\dots,v_{m}$ are extreme points of $P_{1}$.
Now we apply Lemma \ref{gincap} to $P_1$. 
Let $\eta_2 < \min\{\varepsilon_1,  \gamma(P_1, \delta_1, \eta_1)\}$.
By Lemma \ref{gincap} there exists $z_2 \in P_1$, $\|v_{2}-z_{2}\|\leq\eta_2$,  and $0 < r_2 \leq \eta_2$ such that $B^n_2(z_2,r_2) \subset P_1$
and such that 
$$
Q_2 = \conv\left(\mathcal{P}_{\varepsilon_1}, B^n_2(z_2,r_2), v_3, \dots , v_m\right) 
$$
has  $v_3, \dots , v_m$ as  extreme points, 
\begin{equation}\label{dq2p1}
d_H(Q_2,P_1) \leq \eta_2,
\end{equation} 
and for sufficiently small $\delta_2$, 
\begin{equation}\label{v2-q}
\|v_2- g\left(Q_2 \setminus (Q_2)_{\delta_2 } \right)\| \leq 2\eta_2.
\end{equation}  
As $\|v_{1}-z_{1}\|\leq\eta_1 $ and $\|v_{2}-z_{2}\|\leq\eta_2 $,  we have that $d_H(Q_2,P) \leq \eta_1 $.
Moreover, as $d_H(Q_2,P_1) \leq \eta_2  < \gamma(P_1, \delta_1, \eta_1)$,  we get by  (\ref{cor1}) with $\varepsilon=\eta_1$ and  by (\ref{v1p1}) 
 that 
$$
\left\| v_1- g(Q_2\setminus (Q_2)_{\delta_1})\right\| \leq \left\| v_1- g(P_1\setminus (P_1)_{\delta_1})\right\| +
\left\| g(P_1\setminus (P_1)_{\delta_1})- g(Q_2\setminus (Q_2)_{\delta_1})\right\| \leq 4 \eta_1.
$$
Now we let $\varepsilon_2< \min\{ \eta_2, \gamma(Q_2, \delta_1, \eta_1)\}$, choose an $\varepsilon_2$-net $\mathcal{P}_{\varepsilon_2}$ on $\partial \left(B^n_2(z_2,r_2) \right)$ and put 
$$
P_2=\conv\left( \mathcal{P}_{\varepsilon_1}, \mathcal{P}_{\varepsilon_2}, v_3, \dots, v_m\right).
$$
Then $P_2 \subseteq Q_2 \subseteq P$ and  $d_H(P_2,Q_2) \leq \varepsilon_2 $.
By (\ref{cor1}),  with $\varepsilon=\eta_2$, and  if we choose in addition $\varepsilon_2 < \eta(Q_2, \delta_2, \eta_2)$, we get
$$\left\| g(Q_2\setminus (Q_2)_{\delta_2}) - g(P_2\setminus (P_2)_{\delta_2})\right\| < \eta_2$$
and  thus, together with (\ref{v2-q}),
\begin{equation}\label{v2p2}
\left\| v_2- g(P_2\setminus (P_2)_{\delta_2})\right\| \leq 3\ \eta_2.
\end{equation}
Please note that $v_{3},\dots,v_{m}$ are extreme points of $P_{2}$.
Now we apply Lemma \ref{gincap} to $P_2$. 
Let $\eta_3 < \min\{\varepsilon_2,  \gamma(P_2, \delta_2, \eta_2)\}$.
By Lemma \ref{gincap} there exists $z_3 \in P_2$, $\|v_{3}-z_{3}\|\leq\eta_3$,  and $0 < r_3 \leq \eta_3$ such that $B^n_2(z_3,r_3) \subset P_2$
and such that 
$$
Q_3 = \conv\left(\mathcal{P}_{\varepsilon_1}, \mathcal{P}_{\varepsilon_1}, B^n_2(z_3,r_3), v_4, \dots , v_m\right) 
$$
has  $v_4, \dots , v_m$ as  extreme points, 
\begin{equation}\label{dq2p1}
d_H(Q_3,P_2) \leq \eta_3,
\end{equation} 
and for sufficiently small $\delta_3$, 
\begin{equation}\label{v3-q}
\|v_3- g\left(Q_3 \setminus (Q_3)_{\delta_3 } \right)\| \leq 2  \  \eta_3.
\end{equation}  
As $\|v_{1}-z_{1}\|\leq\eta_1 $, $\|v_{2}-z_{2}\|\leq\eta_2 $ and $\|v_{3}-z_{3}\|\leq\eta_3 $ we have that $d_H(Q_3,P) \leq \eta_1 $.
Moreover, as $d_H(Q_3,P_2) \leq \eta_3  < \gamma(P_2, \delta_2, \eta_2)$,  we get by  (\ref{cor1}) with $\varepsilon=\eta_2$ and (\ref{v2p2}) 
 that 
$$
\left\| v_2- g(Q_3\setminus (Q_3)_{\delta_2})\right\| \leq \left\| v_2- g(P_2\setminus (P_2)_{\delta_2})\right\| +
\left\| g(P_2\setminus (P_2)_{\delta_2})- g(Q_3\setminus (Q_3)_{\delta_2})\right\| \leq 4 \eta_2.
$$
As $d_H(Q_2,Q_3) \leq \varepsilon_2 < \gamma(Q_2, \delta_1, \eta_1)$, it follows from  (\ref{cor1}) with $\varepsilon=\eta_1$ that $$\left\| g(Q_2\setminus (Q_2)_{\delta_1})- g(Q_3\setminus (Q_3)_{\delta_1})\right\| \leq \eta_1.$$ 
By (\ref{dq2p1}), it also follows from (\ref{cor1}) with $\varepsilon=\eta_1$ that 
$$
 \left\| g(P_1\setminus (P_1)_{\delta_1}) - g(Q_2\setminus (Q_2)_{\delta_1})\right\| \leq \eta_1.
$$
This, together  with (\ref{v1p1}) gives
\begin{eqnarray*}
\left\| v_1- g(Q_3\setminus (Q_3)_{\delta_1})\right\| &\leq& \left\| v_1- g(P_1\setminus (P_1)_{\delta_1})\right\| +
 \left\| g(P_1\setminus (P_1)_{\delta_1}) - g(Q_2\setminus (Q_2)_{\delta_1})\right\| \\
 &&+
\left\| g(Q_2\setminus (Q_2)_{\delta_1})- g(Q_3\setminus (Q_3)_{\delta_1})\right\| \\
&\leq& 5 \eta_1.
\end{eqnarray*}
We continue to obtain  $Q=Q_{n+1}$ and affine invariant points $\Delta_i = g(Q \setminus Q_{\delta_i})$, $1 \leq i \leq n+1$,  such that
for all $i$, 
  $$ \|v_i - \Delta_i(Q)\| \leq (n+2) \eta_1<\eta.
$$  
As for $1 \leq i \leq n+1$, the $v_i$ are  affinely independant,  so are the $\Delta_i$.  
\par
It remains to show that ${\cal O}_n=\{K \in \mathcal{K}_n: \mathfrak{P}_n(K)= \mathbb{R}^n\}$ is open in $(\mathcal{K}_n, d_H)$. Observe that
$K\in {\cal O}_n$ if and only if for some $p_1\dots,p_{n+1}\in {\mathfrak P}_n$ (depending on $K$), $$\vol_n(\conv\big(p_1(K)\dots,p_{n+1}(K)\big)>0. $$
 Since
$L\to \vol(\conv\big(p_1(L)\dots,p_{n+1}(L)\big)$ is continuous on 
${\mathcal K}_n$, it follows that ${\cal O}_n$ is open. $\square$
\vskip 2mm
\noindent
\begin{corollary} For every  $k \in \mathbb{N}$, $0 \leq k \leq n$, there exists a convex body $Q_k$
such that $\mathfrak{P}(Q_k)$ is a $k$-dimensional affine subspace of $\mathbb{R}^n$.
\end{corollary}
\noindent
{\bf Proof.}
For $k=0$, we take a centrally symmetric body. For $k=n$, we take the body $Q$ of Theorem \ref{dicht}. For $1\le k\le n-1$, we take the intermediate bodies $Q_k$ constructed in the proof of Theorem \ref{dicht}.  $\square$

\vskip 4mm

Mathieu Meyer \\
  {\small        Universit\'{e} de Paris Est - Marne-la-Vall\'{e}e}\\
    {\small      Equipe d'Analyse et de Math\'ematiques Appliqu\'{e}es}\\
   {\small       Cit\'e  Descartes - 5, bd Descartes }\\
    {\small      Champs-sur-Marne
          77454 Marne-la-Vall\'{e}e,  France} \\
          {\small \tt mathieu.meyer@univ-mlv.fr} \\

   Carsten Sch\"utt \\
     {\small       Christian Albrechts Universit\"at }\\
        {\small    Mathematisches Seminar }\\
        {\small    24098 Kiel, Germany} \\
         {\small \tt schuett@math.uni-kiel.de}   \\

     \and Elisabeth Werner\\
{\small Department of Mathematics \ \ \ \ \ \ \ \ \ \ \ \ \ \ \ \ \ \ \ Universit\'{e} de Lille 1}\\
{\small Case Western Reserve University \ \ \ \ \ \ \ \ \ \ \ \ \ UFR de Math\'{e}matiques }\\
{\small Cleveland, Ohio 44106, U. S. A. \ \ \ \ \ \ \ \ \ \ \ \ \ \ \ 59655 Villeneuve d'Ascq, France}\\
{\small \tt elisabeth.werner@case.edu}\\ \\


\vskip 5mm

\begin{thebibliography}{~~}
\small



\bibitem{BlaschkeBook}
{\sc W. Blaschke}, {\em Vorlesungen {\"u}ber Differentialgeometrie
II, Affine Differentialgeometrie}. Springer Verlag, Berlin,
(1923).


\bibitem{Boroczky2010}
{\sc K.J. B\"or\"oczky},
{\em Stability of the Blaschke-Santal\'o and the affine isoperimetric inequality}
Adv. in Math. {\bf 225} (2010), 1914--1928.


\bibitem{BLYZ2012}
{\sc K.J. B\"or\"oczky, E. Lutwak, D. Yang, and G. Zhang}, {\em The logarithmic Minkowski problem}, 
to appear in Journal of AMS.

\bibitem{BourgainMilman1987}
{\sc J. Bourgain and  V. D. Milman}, {\em New volume ratio properties for
convex symmetric bodies in $\mathbb R^n$ }, Invent.\ Math. {\bf 88} (1987),
319--340.


\bibitem{CampiGronchi} {\sc S. Campi and P. Gronchi}, {\em The $L^p$-Busemann-Petty
centroid inequality}, Adv. in Math. {\bf 167} (2002), 128--141.


\bibitem{Ga1}
{\sc  R. J. Gardner}, { \em A positive answer to the Busemann-Petty problem in three dimensions}, 
Ann. of Math. (2) {\bf 140} (1994), 435-47.  

\bibitem{GardnerBook}
{\sc R.J. Gardner}
{\em Geometric tomography.}
 Second edition. Encyclopedia of Mathematics and its Applications, 58. Cambridge University Press, Cambridge (2006).

\bibitem{Ga3}
{\sc R. J. Gardner}, {\em The dual Brunn-Minkowski theory for bounded Borel sets: Dual affine 
quermassintegrals and inequalities}, Adv. Math. {\bf 216} (2007), 358-386. 



\bibitem{GaKoSch}
{\sc  R. J. Gardner, A. Koldobsky, and T.  Schlumprecht}, {\em An analytical solution 
to the Busemann-Petty problem on sections of convex bodies}, Ann. of Math. (2) { \bf 149} (1999), 691-703.



\bibitem{GaZ}
{\sc R. J. Gardner and G. Zhang},
{\em Affine inequalities and radial mean bodies.}
 Amer. J. Math. {\bf 120}, no.3, (1998), 505-528.

\bibitem{GMR1988}
{\sc Y. Gordon, M. Meyer and S. Reisner}, {\em Zonoids with minimal
volume product--a new proof}, Proc.\ Amer.\ Math.\ Soc.  {\bf 104} (1988),
273--276.

\bibitem{GrZh}
{\sc
E. Grinberg and G. Zhang,} {\em Convolutions, transforms, and convex bodies}, Proc. 
London Math. Soc. (3) {\bf 78} (1999), 77--115.




\bibitem{Gruenbaum1963}
{\sc B. Gr\"unbaum},
{\em Measures of symmetry for convex sets},
Proc. Sympos. Pure Math. 7, (1963), 233--270.


\bibitem{Hab}
{\sc C. Haberl}, {\em Blaschke valuations}, Amer. J.  Math., 133, (2011), 717--751.




\bibitem{HabSch2}
{\sc C. Haberl and F.  Schuster,} {\em General Lp affine isoperimetric inequalities}.
J. Differential Geometry {\bf 83} (2009), 1-26.

\bibitem{HLYZ}
{\sc C. Haberl, E. Lutwak, D. Yang and G. Zhang,} {\em The even Orlicz Minkowski problem}, 
Adv. Math. {\bf 224} (2010), 2485-2510.

\bibitem{Klain1}
{\sc D. Klain},
{\em  Star valuations and dual mixed volumes}, Adv. Math. {\bf 121} (1996), 80-101. 

\bibitem{Klain2}
{\sc D. Klain},
{\em Invariant valuations on star-shaped sets},   Adv. Math.  {\bf 125} (1997), 95-113. 




\bibitem{Kuperberg2008}
{\sc G. Kuperberg}, {\em  From the Mahler Conjecture to Gau{\ss}
Linking Integrals}, Geometric And Functional Analysis {\bf 18} (2008), 870--892.



\bibitem{Lud2}
{\sc M. Ludwig}, {\em Ellipsoids and matrix valued valuations}, Duke Math. J. {\bf 119} (2003), 159-188.


\bibitem{Lud3}
{\sc M. Ludwig}, {\em
Minkowski areas and valuations}, 
J. Differential Geometry, 86, (2010), 133--162.


\bibitem{LR1}
{\sc M. Ludwig and M. Reitzner,}  {\em A Characterization of Affine Surface Area}, Adv. in Math. {\bf 147} (1999), 138-172.

\bibitem{LR2}
{\sc M. Ludwig and M. Reitzner,} {\em A classification of $SL(n)$
invariant valuations.}  Annals of Math. {\bf 172 } (2010), 1223-1271. 



\bibitem{Lu1}{\sc E. Lutwak}, {\em The Brunn-Minkowski-Firey theory I :  Mixed volumes and the Minkowski problem},
J. Differential Geom. {\bf 38} (1993), 131--150.

\bibitem{Lu2}{\sc E. Lutwak}, {\em The Brunn-Minkowski-Firey theory II : Affine and
geominimal surface areas}, Adv. Math. {\bf 118}  (1996),   244-294.

\bibitem{LutwakZhang1997} {\sc E. Lutwak and G. Zhang}, {\em Blaschke-Santal\'{o}
inequalities}, J. Differential Geom.  47,  (1997), 1--16.


\bibitem{LutwakYangZhang200/1} {\sc E. Lutwak, D. Yang and G. Zhang}, {\em $L^p$ affine isoperimetric
inequalities}, J. Differential Geom. {\bf 56} (2000), 111--132.



\bibitem{LYZ2000} {\sc E. Lutwak, D. Yang and G. Zhang}, {\em A new ellipsoid associated with convex bodies}, Duke Math. J. { \bf 104} (2000), 375--390.

\bibitem{LYZ2002} {\sc E. Lutwak, D. Yang and G. Zhang}, {\em Sharp Affine $L_p$ Sobolev inequalities}, 
J. Differential Geom. {\bf 62} (2002), 17--38.


\bibitem{LYZ2002/1} {\sc E. Lutwak, D. Yang and G. Zhang}, {\em The Cramer--Rao inequality for star bodies}, Duke Math. J. {\bf 112} (2002), 59-81.

\bibitem{LYZ2004} {\sc E. Lutwak, D. Yang and G. Zhang}, {\em Volume inequalities for subspaces of $L_p$},  
J. Differential Geom. {\bf 68} (2004), 159--184.



\bibitem{MeyerWerner1998} {\sc M. Meyer and E. Werner}, {\em The Santalo-regions of a convex body},
Transactions of the AMS {\bf 350}   (1998), 4569--4591.


\bibitem{MeyerWerner2000} 
{\sc M. Meyer and E. Werner},
{\em On the \ensuremath{p}-affine surface area.}
Adv. Math. {\bf 152} (2000), 288--313.



\bibitem{MeyerSchuettWerner2011} 
{\sc M. Meyer, C. Sch\"utt and E. Werner},
{\em New affine measures of symmetry for convex bodies}, 
Adv. Math. {\bf 228}, (2011), 2920--2942.




\bibitem{Nazarov}
{\sc F. Nazarov}, {\em The H\"ormander proof of the Bourgain-Milman theorem}, preprint, 2009.


\bibitem{NPRZ2010}
{\sc F. Nazarov,  F. Petrov, D. Ryabogin and A. Zvavitch},  {\em A remark on the Mahler conjecture: local minimality of the unit cube}, 
Duke Mathematical J. {\bf 154}, (2010), 419--430.


\bibitem{Paouris2006} {\sc G. Paouris},
{\em Concentration of mass on convex bodies}, Geometric and
Functional Analysis {\bf 16},  (2006), 1021--1049.

\bibitem{Petty1972}
{\sc C. M. Petty}, {\em Isoperimetric problems}, Proc. Conf. Convexty and Combinatorial
Geometry Univ. Oklahoma 1971, University of Oklahoma, (1972),  26--41.






\bibitem{R1}
{\sc S. Reisner}, {\em Zonoids with minimal volume-product}, Math.\ Z. {\bf 192} (1986),
339--346.

\bibitem{R2}
{\sc S. Reisner},  {\em Minimal volume product in Banach spaces with
a 1-unconditional basis}, J. London Math.\ Soc.{\bf 36} (1987), 126--136.

\bibitem{ReisnerSchuettWerner2011}
{\sc S. Reisner, C. Sch\"utt and E. Werner}, 
{\em Mahler's conjecture and curvature} 
International Mathematics Research
Notices, DOI 10.1093/imrn/rnr003 (2011)




\bibitem{S-R}
{\sc J. Saint-Raymond}, {\em Sur le volume des corps convexes
sym\'etriques},  S\'eminaire d'Initiation \`a l'Analyse, 1980-1981, Universit\'e
PARIS VI, Paris 1981.



\bibitem{SchneiderBuch} {\sc R. Schneider}, 
{\em Convex Bodies: The
Brunn-Minkowski Theory}, Encyclopedia of Mathematics and its
Applications  44, Cambridge University Press, Cambridge (1993).



\bibitem{Schuster2010}
{\sc
F. Schuster}, {\em Crofton measures and Minkowski valuations}, Duke Math. J. {\bf 154} (2010), 1--30.

\bibitem{SchusterWeberndorfer2012}
{\sc
F. Schuster and M. Weberndorfer}, {\em
Volume Inequalities for Asymmetric Wulff Shapes},  
J. Differential Geom., in press.


\bibitem{SchuettWerner1990} {\sc C. Sch\"utt and E. Werner}, {\em The convex floating body},   
{Math. Scand.} {\bf 66}, (1990), 275--290.

\bibitem{SchuettWerner1992} {\sc C. Sch\"utt and E. Werner}
{\em The convex floating body of almost polygonal bodies}, 
Geom. Dedic. {\bf 44}, (1992), 169--188.


\bibitem{SchuettWerner2003} 
{\sc C. Sch\"utt and E. Werner}, 
{\em Polytopes with Vertices Chosen Randomly
from the Boundary of a Convex Body},   Geom. Aspects of Funct. Analysis,
Lecture Notes in Math. {\bf 1807},  (2003), 241--422.




\bibitem{SW5}
{\sc C. Sch{\"u}tt and E. Werner}, {\em Surface bodies and
p-affine surface area.} Adv. Math. {\bf 187}  (2004), 98-145.



\bibitem{spingarn}
{\sc J. E. Spingarn},
{\em An Inequality for Sections and Projections of a Convex Set}, Proc.\ Amer.\ Math.\ Soc., {\bf 118}, No. 4, (1993), 1219--1224.



\bibitem{SA1}
{\sc A. Stancu,} {\em The Discrete Planar $L_0$-Minkowski
Problem.} Adv. Math. {\bf 167} (2002),  160-174.

\bibitem{SA2}
{\sc A. Stancu}, {\em On the number of solutions to the
discrete two-dimensional $L_0$-Minkowski problem.} Adv. Math. {\bf
180} (2003), 290-323.



\bibitem{WernerYe2008} {\sc E. Werner and D. Ye}, {\em New $L_{p}$ affine isoperimetric inequalities},
Adv. Math. {\bf 218} (2008), no. 3, 762-780.


\bibitem{WernerYe2010} {\sc E. Werner and  D. Ye}, {\em Inequalities for mixed $p$-affine surface area},
{Math. Ann.} {\bf  347}  (2010), 703-737.

\bibitem{Zhang1991}
{\sc G. Zhang}, {\em
Restricted chord projection and affine inequalities}, Geom. Dedicata,
{\bf 39} (1991), 213--222.

\bibitem{Z1}{\sc G. Zhang}, {\em Intersection bodies and Busemann-Petty inequalities in $\mathbb{R}^4$}, Annals of Math. {\bf 140} (1994), 331-346. 


\bibitem{Z2}{\sc G. Zhang}, {\em A positive answer to the Busemann-Petty problem in four dimensions}, Annals of. Math. {\bf 149} (1999), 535-543.

\bibitem{Z3}{\sc G. Zhang}, {\em New Affine  Isoperimetric Inequalities}, ICCM 2007, Vol. II, 239-267.







\end{thebibliography}
\end{document}